\documentclass[preprint,11pt]{elsarticle}
\usepackage[top=1.3in, bottom=2cm, left=2.5cm, right=2.5cm]{geometry}

\usepackage{graphics,epsfig}

\usepackage{amssymb,amsmath}
\usepackage{graphicx}

\journal{}
\begin{document}
\begin{frontmatter}
\newtheorem{thm}{Theorem}[section]
\newtheorem{lem}{Lemma}[section]
\newtheorem{coro}{Corollary}[section]
\newtheorem{prop}{Proposition}[section]
\newdefinition{definition}{Definition}
\newdefinition{rmk}{Remark}
\newproof{pf}{Proof}
\newtheorem{ex}{Example}
\title{A Sharp Upper Bound for Region Unknotting Number of Torus Knots\tnoteref{t1}}
\tnotetext[t1]{This document is a collaborative effort.}

\author[v]{Vikash. S}

\author[mp]{Madeti. P\corref{cor1}}
\ead{prabhakar@iitrpr.ac.in}
\ead[url]{http://www.iitrpr.ac.in/html/faculty/prabhakar.shtml}
\cortext[cor1]{Corresponding author}
\address[v]{Department of Mathematics, IIT Ropar, Rupnagar- 140001, India.}
\address[mp]{Room No. 207, Department of Mathematics, IIT Ropar, Rupnagar - 140001, India.}

\begin{abstract}
Region crossing change for a knot or a proper link is an unknotting operation. In this paper, we provide a sharp upper bound on the region unknotting number for a large class of torus knots and proper links. Also, we discuss conditions on torus links to be proper. 
\end{abstract}

\end{frontmatter}



\section{Introduction} The unknotting number $u(D)$ of a knot diagram $D$ is the minimal number of crossing changes required to unknot $D$ and the unknotting number $u(K)$ of a knot $K$ is given by $u(K) = inf\{u(D) : D\ is\ a\ knot\ diagram\ of\ K\}$. It is not easy to calculate $u(K)$, because it is taken over an infinite number of diagrams of $K$. Unknotting number lead to consider different types of unknotting operations i.e., the local transformations on link diagrams which can transform the knot diagram into a trivial knot diagram. It was shown that $\sharp$-operation~\cite{hash}, $\delta$-operation~\cite{delta}, $3$-gon operation~\cite{3-gon}, H(n)-operation~\cite{H(n)} and $n$-gon~\cite{n-gon} operations are unknotting operations.

Recently, Ayaka Shimizu~\cite{ayaka} proposed an unknotting operation for knots called \textit{region crossing change}. In this operation, $K$ to be a knot and $D$ be any diagram of $K$, then a region crossing change at a region $R$ of diagram $D$ is defined to be the crossing changes at all the crossing points on $\partial R$. The region unknotting number $u_R (D)$ of a knot diagram $D$ is the minimal number of region crossing changes required to transform $D$ into  a diagram of the trivial knot without Reidemeister moves. The region unknotting number $u_R(K)$ of $K$ is defined to be the minimal $u_R(D)$ taken over all minimal diagrams $D$ of $K$. 

In general, for links, the region crossing change is not an unknotting operation. In \cite{ayaka}, Ayaka Shimizu considered the standrard diagram of Hopf link and showed that it cannot be transformed into a diagram of trivial link by region crossing changes. In \cite{incidence} and \cite{proper} Cheng Zhiyun defined a link $L= K_1 \bigcup K_2 \bigcup \cdots \bigcup K_n$ to be proper if for each $i,\ 1\leq i \leq n$, $\sum_{j\neq i}lk(K_i,K_j) \equiv 0\ (mod\ 2)$  and proved that \textit{region crossing change on a link is an unknotting operation if and only if the link is proper}.

The region unknotting number is intuitive and attractive for knot theorists and it is not surprising that there is no algorithm to compute it for a given knot or proper link at the moment. In \cite{ayaka}, \cite{ayaka2}, Ayaka Shimizu showed that for a twist knot $K$, $u_R(K)=1$ and for torus knots of type $K(2,4m\pm1)$, $u_R(K(2,4m\pm1))=m$, where $m\in \mathbb{Z^+}$. As the unknotting number for torus knots is already known \cite{km-2}, \cite{km-3}, it would be worth attempting to find the region unknotting number for torus knots. In this paper, we discuss the condition for torus links to be proper and provide a sharp upper bound for a large class of torus knots and proper torus links, including all $2-, 3-, 4-, 5-, 6-$ and $7-$braid torus knots/links.

In Section \ref{properlinks}, we discuss the condition for torus links to be proper. In Section \ref{pre}, we discuss basic results, depending on elementary braid relations, which helps to provide a sharp upper bound for region unknotting number for torus links. In Section \ref{ubound}, we provide a sharp upper bound for region unknotting number of $K(p,np+a)$ torus knots/links, where $a=1,2,3,4, -1, -2$,  and partially for $K(p,np+5)$ type of torus links.
\section{Conditions for torus link to be proper}\label{properlinks}
In \cite{proper}, Cheng Zhiyun showed that, \textit{region crossing change is an unknotting operation on a link $L= K_1 \bigcup K_2 \bigcup \cdots \bigcup K_n$ if and only if $\sum_{j\neq i}lk(K_i,K_j) \equiv 0\ (mod\ 2)$ for each $1\leq i \leq n$}. If $gcd(p,q)=d$, we use the notation $K(p,q) = K_1 \bigcup K_2 \bigcup \cdots \bigcup K_d$, where each $K_i$ is a component of $K(p,q)$. To provide the condition for torus links to be proper, we consider the diagram for  $K(p,q)$ as a closure of toric braid $(\sigma_1\sigma_2\cdots\sigma_{p-1})^q$. Also, we denote $\displaystyle \sum_{i=1}^n i$ by $\sum n$. Given a $d$-component torus link $K(p,q)$, the number of crossings between any two components is same i.e., for $i\neq j$, $l\neq m$, $\sharp
(K_i\bigcap K_j) = \sharp (K_l\bigcap K_m)$. Also, it is worth noting that, for $i\neq j$, $lk(K_i,K_j) = \frac{1}{2}\sharp(K_i\bigcap K_j)$.
%
\begin{thm}\label{thm1} 
A torus link $K(p,q)$ is proper iff \[\frac{pq(d-1)}{d^2}\equiv 0(mod\ 2),\] where: d=gcd(p,q).
\end{thm}
\begin{pf}
Consider a $d$-component torus link $K(p,q)$. Since we are considering $K(p,q)$ as  a closure of toric braid 
$(\sigma_1\sigma_2\cdots\sigma_{p-1})^q$, each component $K_i$ contains $p/d$ strands. For each bracket $(\sigma_1\sigma_2\cdots\sigma_{p-1})$, there exists only one $i$ such that $\sharp (K_i \bigcap K_i) = (\frac{p}{d}-1)$. 
Hence the total number of self-intersections is
\[\displaystyle \sum_{i=1}^{d}\sharp \left(K_i \bigcap K_i\right)  = q\left(\frac{p}{d}-1\right).\]
Since,
\begin{align}\displaystyle\sum_{i\neq j}  \sharp(K_i \bigcap K_j)  &= total\ no\ of\ crossings - \displaystyle \sum_{i=1}^{d} \sharp(K_i \bigcap K_i), \notag 
\end{align}
we have 
\begin{align}  \displaystyle\sum_{i\neq j}  \sharp(K_i \bigcap K_j)=\frac{pq(d-1)}{d}.
\end{align}
Because $\sharp(K_i \bigcap K_j) = \sharp(K_l \bigcap K_m)$, for any $i\neq j, l\neq m$, we have
\begin{align}\displaystyle \sum_{i\neq j} \sharp(K_i \bigcap K_j) = \frac{(d-1)d}{2}\left(\sharp(K_1 \bigcap K_2)\right)\tag{2}.
\end{align}
Thus, from (1) and (2), for any $i\neq j$ \[\sharp(K_i \bigcap K_j)=\frac{2pq}{d^2}.\]
Since all crossings in $K(p,q)$ are positive, for fixed $n$
\begin{align}\displaystyle \sum_{j\neq n} lk(K_n,K_j) &= \displaystyle \frac{1}{2}\sum_{j\neq n}\sharp  (K_n \bigcap K_j)\notag \\ 
&=\frac{1}{2}(d-1)\left(\sharp (K_i \bigcap K_j)\right) \notag \\
&=\frac{pq(d-1)}{d^2}. \notag
\end{align}
Hence, a torus link $K(p,q)$, where d=gcd(p,q), is proper iff\\
$~~~~~~~~~~~~~~~~~~~~~~~~~~~~~~~~~~~~~~~~~~~~~~~~\displaystyle \frac{pq(d-1)}{d^2}\equiv 0(mod\ 2).$ \hfill $\square$
\end{pf}
\begin{coro}\label{coro1}
Let $K(p,q)$ be any d-componant torus link where $p=2^mk$ and $q=2^nk'$, where $m$ and $n$ are the highest powers of $2$ in $p$ and $q$, respectively. Then $K(p,q)$ is proper iff either $m=n=0$ or $m\neq n$.
\end{coro}
\section{Preliminaries}\label{pre}
To find the sharp upper bounds on the region unknotting number of torus links, we use the following:
\begin{thm}\label{thm2} (\cite{ps1}, Theorem 3.1)
 For every $p$, the $p$-braid
\[\underbrace{\sigma_1\sigma_2\cdots\sigma_{p-1}}_1\underbrace{\sigma_1\sigma_2\cdots\sigma_{p-2}\sigma_{p-1}^{-1}}_2 
\underbrace{\sigma_1\sigma_2\cdots\sigma_{p-2}^{-1}\sigma_{p-1}^{-1}}_3\cdots \underbrace{\sigma_1^{-1}\sigma_2^{-2}\cdots\sigma_{p-1}^{-1}}_p\]
is a trivial $p$-braid.
\end{thm}
\begin{rmk}\label{coro2} Using Theorem \ref{thm2}, we observe the following:\\
(i) For every $p$, the $p$-braid
\[\underbrace{\sigma_1^{-1}\sigma_2^{-2}\cdots\sigma_{p-1}^{-1}}_1\underbrace{\sigma_1^{-1}\sigma_2^{-2}\cdots\sigma_{p-2}^{-1}
\sigma_{p-1}}_2\cdots\underbrace{\sigma_1^{-1}\sigma_2\cdots\sigma_{p-1}}_{p-1}\underbrace{\sigma_1\sigma_2\cdots
\sigma_{p-1}}_p \]
is a trivial $p$-braid.\\
(ii) The closure of $n+1$ braids 
\[\underbrace{\sigma_1\sigma_2\cdots\sigma_{n-1}\sigma_{n}^{-1}}_1\underbrace{\sigma_1\sigma_2\cdots\sigma_{n-1}^{-1}
\sigma_{n}^{-1}}_2\cdots \underbrace{\sigma_1^{-1}\sigma_2^{-1}\cdots\sigma_{n}^{-1}}_n \]
and
\[\underbrace{\sigma_1\sigma_2\cdots\sigma_{n-1}\sigma_{n}}_1\underbrace{\sigma_1\sigma_2\cdots\sigma_{n-1}
\sigma_{n}^{-1}}_2\cdots \underbrace{\sigma_1\sigma_2^{-1}\cdots\sigma_{n}^{-1}}_n \]
are trivial links.
\end{rmk}
\begin{lem}\label{thm7}  For every $p$ and $a$,\ where\ $p>a$, the\ $p$-braid\\
$\underbrace{\eta_1\kappa_{p-1}}_1\underbrace{\eta_2\kappa_{p-2}\sigma_{p-1}^{-1}}_2\cdots \underbrace{\eta_{a-1}\kappa_{p-a+1}\sigma_{{p-a}+2}^{-1}\cdots\sigma_{p-1}^{-1}}_{a-1}\underbrace{\eta_a\sigma_{{p-a}+1}^{-1}\cdots\sigma_{p-1}^{-1}}_a$ is Markov equivalent of $\eta_1'\eta_2'\cdots\eta_a'$, where
$\eta_i'=\sigma_1^{g_{i,1}}\sigma_2^{g_{i,2}}\cdots \sigma_{{p-a}-1}^{g_{i,{p-a}-1}}$, $\eta_i = \eta_i'\sigma_{{p-a}}^{g_{i,{p-a}}}$ with $g_{i,j}= \pm 1$, for $i=1,2,\ldots,a$; $j=1,2,\ldots,{p-a}$ and $\kappa_{j}=\sigma_{p-a+1}\sigma_{{p-a}+2}\cdots \sigma_j$.
\end{lem}
\begin{pf} Using fundamental braid relations, we observe that
\begin{align}
&\underbrace{\eta_1\kappa_{p-1}}_1\underbrace{\eta_2\kappa_{p-2}\sigma_{p-1}^{-1}}_2\cdots \underbrace{\eta_{a-1}\kappa_{p-a+1}\sigma_{{p-a}+2}^{-1}\cdots\sigma_{p-1}^{-1}}_{a-1}\underbrace{\eta_a\sigma_{{p-a}+1}^{-1}\cdots\sigma_{p-1}^{-1}}_a \notag \\
&\sim_M \underbrace{\eta_1'\sigma_{p-a}^{g_{1,p-a}}\kappa_{p-1}}_1\underbrace{\eta_2'\sigma_{p-a}^{g_{2,p-a}}\kappa_{p-2}\sigma_{p-1}^{-1}}_2\cdots \underbrace{\eta_a'\sigma_{p-a}^{g_{a,p-a}}\sigma_{{p-a}+1}^{-1}\cdots\sigma_{p-1}^{-1}}_a \notag \\
&\sim_M \underbrace{\eta_1'\sigma_{p-a}^{g_{1,p-a}}\kappa_{p-2}}_1\underbrace{\eta_2'\sigma_{p-a}^{g_{2,p-a}}\kappa_{p-3}\sigma_{p-2}^{-1}}_2\underline{\sigma_{p-1}\sigma_{p-2}}\cdots \underbrace{\eta_a'\sigma_{p-a}^{g_{a,p-a}}\sigma_{{p-a}+1}^{-1}\cdots\sigma_{p-1}^{-1}}_a \notag \\
&~~~~~~~~~~~~~~~~~~~~~~~~~~~~~~~~~~~~~~~~~~~~\vdots \notag \\
&\sim_M \underbrace{\eta_1'\sigma_{p-a}^{g_{1,p-a}}\kappa_{p-2}}_1\underbrace{\eta_2'\sigma_{p-a}^{g_{2,p-a}}\kappa_{p-3}\sigma_{p-2}^{-1}}_2\cdots \underbrace{\eta_a'\sigma_{p-a}^{-1}\cdots\sigma_{p-2}^{-1}}_a \underline{\sigma_{p-1}^{g_{a,p-a}}\sigma_{p-2}\cdots\sigma_{p-a}}\notag \\
&\sim_M \underbrace{\eta_1'\sigma_{p-a}^{g_{1,p-a}}\kappa_{p-2}}_1\underbrace{\eta_2'\sigma_{p-a}^{g_{2,p-a}}\kappa_{p-3}\sigma_{p-2}^{-1}}_2\cdots\underbrace{\eta_{a-1}'\sigma_{p-a}^{g_{a-1,p-a}}\sigma_{p-a+1}^{-1}\sigma_{{p-a}+2}^{-1}\cdots\sigma_{p-2}^{-1}}_{a-1}\eta_a'. \notag
\end{align}
Continuing the same procedure after $a-2$ steps, we obtain
\begin{align}
&\sim_M \underbrace{\eta_1'\sigma_{p-a}^{g_{1,p-a}}\kappa_{p-3}}_1\underbrace{\eta_2'\sigma_{p-a}^{g_{2,p-a}}\kappa_{p-4}\sigma_{p-3}^{-1}}_2\cdots\underbrace{\eta_{a-2}'\sigma_{p-a}^{g_{a-2,p-a}}\sigma_{p-a+1}^{-1}\cdots\sigma_{p-3}^{-1}}_{a-2}\eta_{a-1}'\eta_a'. \notag 
\end{align}
After $\sum (a-3)$ steps, this $p$-braid is Markov equivalent of
$ \eta_1'\eta_2'\cdots\eta_a'.$\hfill $\square$
\end{pf}
\begin{lem}\label{thm10} For any two even numbers $p$ and $q$ with $p>q$,
\[\eta_1\kappa_1\eta_2\kappa_2\cdots\eta_q\kappa_q \sim_M \eta_1'\eta_2'\cdots\eta_q',\]
where $\eta_i = \eta_i'\sigma_{p-q}^{g_{i,p-q}}$ and $\eta_i' = \sigma_1^{g_{i,1}}\sigma_2^{g_{i,2}}\cdots
\sigma_{p-q-1}^{g_{i,p-q-1}};$ and $\kappa_i = \sigma_{p-q+1}^{-1}\sigma_{p-q+2}^{-1}\cdots \sigma_{p-i}^{-1}\sigma_{p-i+1}\cdots\sigma_{p-1}$ with $g_{i,j} = \pm 1$ for $i=1,2,\cdots, q$ and $j=1,2,\cdots, p-q$.
\end{lem}
\begin{pf} Using fundamental braid relations we observe that
\begin{align}
&\underbrace{\eta_1\sigma_{p-q+1}^{-1}\cdots\sigma_{p-1}^{-1}}_1\underbrace{\eta_2\sigma_{p-q+1}^{-1}\cdots\sigma_{p-2}^{-1}\sigma_{p-1}}_2\cdots\underbrace{\eta_{q}\sigma_{p-q+1}\cdots\sigma_{p-1}}_q \notag \\
&\sim_M\underbrace{\eta_1\sigma_{p-q+1}^{-1}\cdots\sigma_{p-2}^{-1}}_1\underbrace{\eta_2\sigma_{p-q+1}^{-1}\cdots\sigma_{p-3}^{-1}\sigma_{p-2}}_2\underline{\sigma_{p-1}^{-1}\sigma_{p-2}^{-1}}\cdots\underbrace{\eta_{q}\sigma_{p-q+1}\cdots\sigma_{p-1}}_q \notag \\
&~~~~~~~~~~~~~~~~~~~~~~~~~~~~~~~~~~~~~~~~~~~~\vdots \notag \\
&\sim_M\underbrace{\eta_1\sigma_{p-q+1}^{-1}\cdots\sigma_{p-2}^{-1}}_1\underbrace{\eta_2\sigma_{p-q+1}^{-1}\cdots\sigma_{p-3}^{-1}\sigma_{p-2}}_2\cdots\underbrace{\eta_{q}'\sigma_{p-q}\cdots\sigma_{p-2}}_q\underline{\sigma_{p-1}^{-1}\sigma_{p-2}^{-1}\cdots\sigma_{p-q}^{-1}} \notag \\
&\sim_M\underbrace{\eta_1\sigma_{p-q+1}^{-1}\cdots\sigma_{p-2}^{-1}}_1\underbrace{\eta_2\sigma_{p-q+1}^{-1}\cdots\sigma_{p-3}^{-1}\sigma_{p-2}}_2\cdots\underbrace{\eta_{q-1}\sigma_{p-q+1}\cdots\sigma_{p-2}}_{q-1}\eta_q'. \notag 
\end{align}
Continuing the same procedure after $q-2$ steps, we obtain\\
$~~~~~~~~~~~~~~~~~~~~\sim_M\underbrace{\eta_1\sigma_{p-q+1}^{-1}\cdots\sigma_{p-3}^{-1}}_1\underbrace{\eta_2\sigma_{p-q+1}^{-1}\cdots\sigma_{p-4}^{-1}\sigma_{p-3}}_2\cdots\underbrace{\eta_{q-2}\sigma_{p-q+1}\cdots\sigma_{p-3}}_{q-1}\eta_{q-1}'\eta_q. $\\
After $\sum (q-3)$ similar steps, this $p$-braid is Markov equivalent of 
$\eta_1'\eta_2'\cdots\eta_q'.$ \hfill  $\square $
\end{pf}
\begin{lem}\label{prop1} Let $\beta_j = \sigma_1^{g_{j,1}}\sigma_2^{g_{j,2}}\cdots \sigma_{p-5}^{g_{j,p-5}}$, where $g_{j,k} = \pm 1$, $1\leq j\leq 2$. Then\\
\[\beta_1\sigma_{p-4}^{g_{1,p-4}}\sigma_{p-3}\sigma_{p-2}\sigma_{p-1}^{-1}\beta_2\sigma_{p-4}^{g_{2,p-4}}\sigma_{p-3}
\sigma_{p-2}^{-1}\sigma_{p-1}^{-1}\sim_M \beta_1 \beta_2.\] 
\end{lem}
\begin{pf} Consider that
\begin{align}
L.H.S &\sim_M \beta_1\sigma_{p-4}^{g_{1,p-4}}\sigma_{p-3}\sigma_{p-2}\beta_2\sigma_{p-4}^{g_{2,p-4}}\sigma_{p-3}\sigma_{p-2}^{-1}\underline{\sigma_{p-1}^{-1}
\sigma_{p-2}^{-1}} \notag \\
&\sim_M \beta_1\sigma_{p-4}^{g_{1,p-4}}\sigma_{p-3}\beta_2 \sigma_{p-4}^{g_{2,p-4}}\sigma_{p-3}^{-1}\underline{\sigma_{p-3}^{-1}\sigma_{p-2}\sigma_{p-3}} \notag \\
&\sim_M \beta_1\sigma_{p-4}^{g_{1,p-4}}\beta_2\underline{\sigma_{p-4}^{-1}\sigma_{p-3}^{g_{2,p-4}}\sigma_{p-4}}\notag \\
&\sim_M \beta_1\beta_2.\notag
\end{align} 
\end{pf}
\begin{lem}\label{thm12}
If~ $p(\geq 4)\equiv 0\ or\ \pm2\ (mod\ 6)$, then the closure of $\sigma_1^{-1}\sigma_2^{-1}\cdots\sigma_{p-1}^{-1} \sigma_1^{g_{2,1}}\sigma_2^{g_{2,2}}
\cdots\sigma_{p-1}^{g_{2,p-1}}\\ \sigma_1^{g_{3,1}} \sigma_2^{g_{3,2}}\cdots\sigma_{p-1}^{g_{3,p-1}}$,
is a trivial link. Here \\
$g_{2,j} = -1$ for $j \in \{p-2,p-3, p-8, p-9, \cdots, p-2-6(\lfloor\frac{p+2}{6}\rfloor-1),p-3-6(\lfloor\frac{p-4}{6}\rfloor)\}$,\\
$g_{3,j} = -1$ for $j \in \{p-3,p-4, p-9, p-10, \cdots, p-3-6(\lfloor\frac{p+2}{6}\rfloor-1),p-4-6(\lfloor\frac{p-4}{6}\rfloor)\}$.
\end{lem}
\begin{pf} By using fundamental braid relations, we obeserve that the $p$-braid $\sigma_1^{-1}\sigma_2^{-1}\cdots\sigma_{p-1}^{-1}\sigma_1^{g_{2,1}}\\ \sigma_2^{g_{2,2}}
\cdots\sigma_{p-1}^{g_{2,p-1}}\sigma_1^{g_{3,1}}\sigma_2^{g_{3,2}}\cdots\sigma_{p-1}^{g_{3,p-1}}$ is Markov equivalent of\\
\textbf{Case 1}. If $p\equiv 0\ (mod\ 6)$,
\begin{align}
&\underbrace{\sigma_1^{-1}\sigma_2^{-1}\cdots\sigma_{p-2}^{-1}}_1
\underbrace{\sigma_1\sigma_2\cdots\sigma_{p-3}^{-1}\sigma_{p-2}}_2
\underbrace{\sigma_1\sigma_2^{-1}\cdots\sigma_{p-3}\sigma_{p-2}}_3
\underline{\sigma_{p-1}^{-1}\sigma_{p-2}^{-1}\sigma_{p-3}^{-1}} \notag \\
&\sim_M \underbrace{\sigma_1^{-1}\sigma_2^{-1}\cdots\sigma_{p-3}^{-1}}_1
\underbrace{\sigma_1\sigma_2\sigma_{3}^{-1}\cdots\sigma_{p-4}\sigma_{p-3}}_2\underline{\sigma_{p-2}^{-1}
\sigma_{p-3}^{-1}}\underbrace{\sigma_1\sigma_2^{-1}\cdots\sigma_{p-4}^{-1}}_3 \notag \\
&\sim_M \underbrace{\sigma_1^{-1}\sigma_2^{-1}\cdots\sigma_{p-4}^{-1}}_1
\underbrace{\sigma_1\sigma_2\cdots\sigma_{p-4}}_2\underbrace{\sigma_1\sigma_2^{-1}\cdots\sigma_{p-9}^{-1}
\sigma_{p-8}\cdots\sigma_{p-5}\sigma_{p-4}^{-1}}_3 \notag \\
&\sim_M \underbrace{\sigma_1^{-1}\sigma_2^{-1}\cdots\sigma_{p-5}^{-1}}_1
\underbrace{\sigma_1\sigma_2\sigma_{3}^{-1}\cdots\sigma_{p-5}}_2\underbrace{\sigma_1\sigma_2^{-1}
\cdots\sigma_{p-7}\sigma_{p-6}\sigma_{p-5}^{-1}}_3 \underline{\sigma_{p-4}\sigma_{p-5}
\sigma_{p-6}^{-1}}\notag \\
&\sim_M \underbrace{\sigma_1^{-1}\sigma_2^{-1}\cdots\sigma_{p-6}^{-1}}_1
\underbrace{\sigma_1\sigma_2\sigma_{3}^{-1}\cdots\sigma_{p-7}\sigma_{p-6}}_2\underline{\sigma_{p-5}\sigma_{p-6}^{-1}}
\underbrace{\sigma_1\sigma_2^{-1}\cdots\sigma_{p-8}\sigma_{p-7}}_3 \notag \\
&\sim_M \underbrace{\sigma_1^{-1}\sigma_2^{-1}\cdots\sigma_{p-7}^{-1}}_1
\underbrace{\sigma_1\sigma_2\cdots\sigma_{p-9}^{-1}\sigma_{p-8}^{-1}
\sigma_{p-7}}_2\underbrace{\sigma_1\sigma_2^{-1}\cdots\sigma_{p-10}^{-1}\sigma_{p-9}^{-1}\sigma_{p-8}\sigma_{p-7}}_3 \notag\\
&~~~~~~~~~~~~~~~~~~~~~~~~~~~~~~~~~~~~~~~~~~~~\vdots \notag \\
&\sim_M \sigma_1^{-1}\sigma_2^{-1}\sigma_{3}^{-1}\sigma_{4}^{-1}\sigma_{5}^{-1}
\sigma_1\sigma_2\sigma_{3}^{-1}\sigma_{4}^{-1}\sigma_{5}\sigma_1\sigma_2^{-1}\sigma_{3}^{-1}\sigma_{4}\sigma_{5}.\notag
\end{align}
\textbf{Case 2}. If $p\equiv 2\ (mod\ 6)$, then
\[\sigma_1^{-1}\sigma_2^{-1}\sigma_{3}^{-1}\sigma_{4}^{-1}\sigma_{5}^{-1}\sigma_{6}^{-1}\sigma_{7}^{-1}
\sigma_1\sigma_2\sigma_{3}\sigma_{4}\sigma_{5}^{-1}\sigma_{6}^{-1}\sigma_{7}
\sigma_1\sigma_2\sigma_{3}\sigma_{4}^{-1}\sigma_{5}^{-1}\sigma_{6}\sigma_{7}.\]
\textbf{Case 3}. If $p\equiv -2\ (mod\ 6)$, then
\[\sigma_1^{-1}\sigma_2^{-1}\sigma_{3}^{-1}\sigma_1^{-1}\sigma_2^{-1}\sigma_{3}\sigma_1^{-1}\sigma_2\sigma_{3}.\]
It is easy to observe that closures of these braids are trivial links. \hfill $\square$
\end{pf}
\begin{lem}\label{prop2} The closure of the braids
\[\sigma_i\sigma_{i+1}\cdots\sigma_j\sigma_i^{-1}\sigma_{i+1}^{-1}\cdots\sigma_{j}^{-1}; if\ i\leq j\]
and
\[\sigma_i\sigma_{i-1}\cdots\sigma_j\sigma_i^{-1}\sigma_{i-1}^{-1}\cdots\sigma_{j}^{-1}; if\ i\geq j\]
are trivial links.
\end{lem}
\begin{pf} The proof directly follows from the elementary braid relations.
\end{pf}
\begin{rmk}\label{coro4} If we take $i = 1$ and $j=n$ in Lemma \ref{prop2}, then the closures of $n+1$-braids
\[\sigma_1\sigma_2\cdots\sigma_n\sigma_1^{-1}\sigma_2^{-1}\cdots\sigma_{n}^{-1}\]
and
\[\sigma_1^{-1}\sigma_2^{-1}\cdots\sigma_{n}^{-1}\sigma_1\sigma_2\cdots\sigma_n\]
are trivial links.
\end{rmk}
\section{Sharp upper bound for region unknotting number of torus knots}\label{ubound}
Let $K$ be any knot/link and $c(K)$ be the crossings number of $k$. Then by~\cite{ayaka}, the trivial upper bound for region unknotting number of knot/link $K$ is $u_R(K) = (c(K)+2)/2$. In case of torus knots/links $K(p,q)$, this trivial upper bound is $u_R(K) = ((p-1)q+2)/2$.

In most cases, to find a sharp upper bound for region unknotting number of torus links $K(p,q)$ with $p<q$, we consider $K(p,q)$ as a closure of $(\sigma_1\sigma_2\cdots\sigma_{p-1})^q$. In particular, we consider the diagram and mark the regions as shown in Figure 1(a). 
 Another diagram we consider for $K(p,q)$ where (q=np+4), is the closure of the braid $(\sigma_1\sigma_2\cdots\sigma_{p-1})^{np}(\sigma_1\sigma_2\cdots\sigma_{3})^{a}\underbrace{\sigma_4\sigma_{3}\cdots\sigma_1}\underbrace{\sigma_{5}\sigma_{4}\cdots \sigma_2}\cdots\\ \cdots \underbrace{\sigma_{p-1}\sigma_{p-2}\cdots\sigma_{p-4}}$ as shown in Figure 1(b), 
  where $c$ and $c'$ are crossing numbers of $K(p,np+4)$ and $K(p,np)$, respectively.

\begin{figure}\label{fig1}
\begin{center}
\includegraphics[width=5cm,height=8cm]{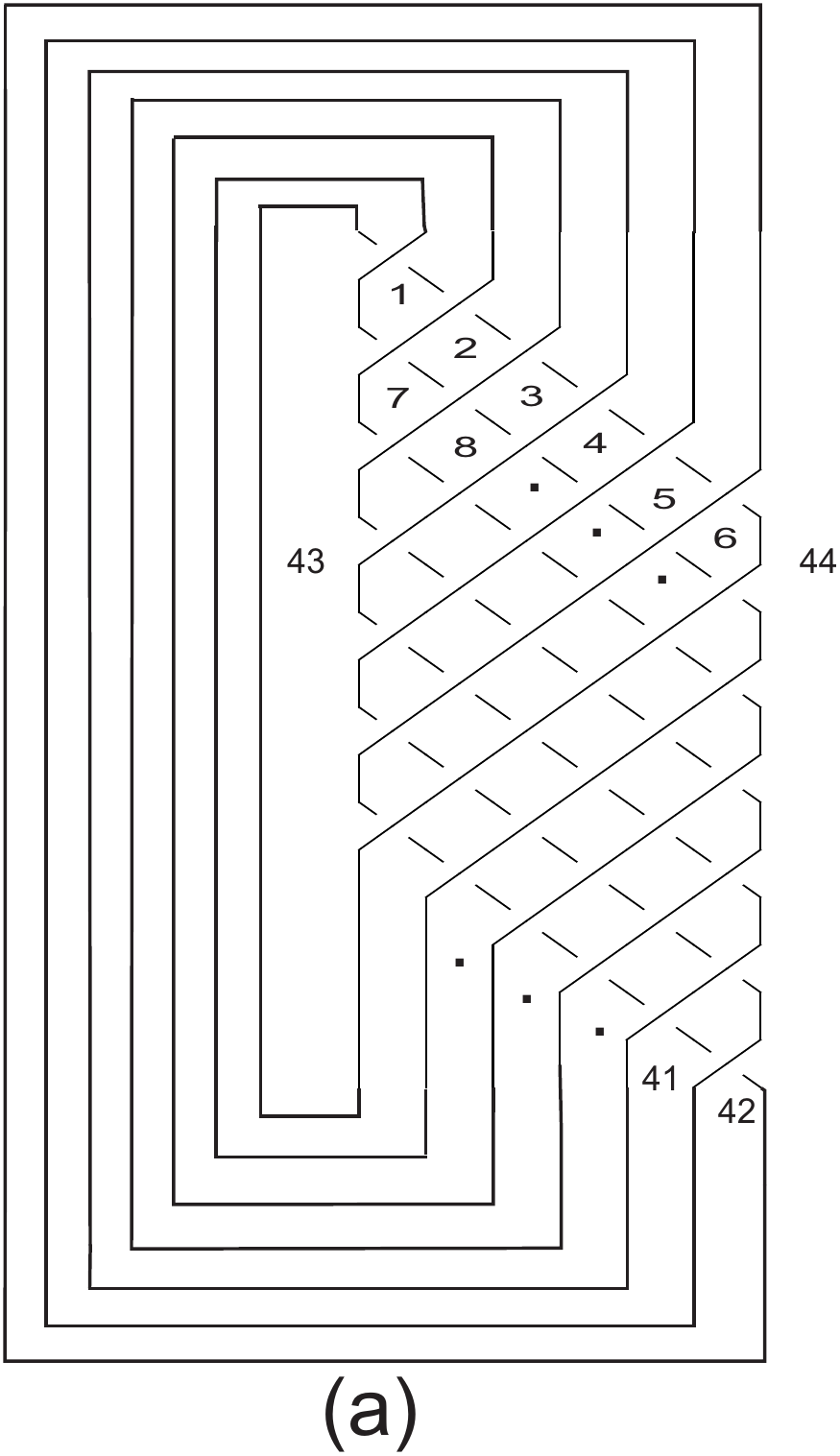}~~~~~~~~~~\includegraphics[width=7cm,height=8cm]{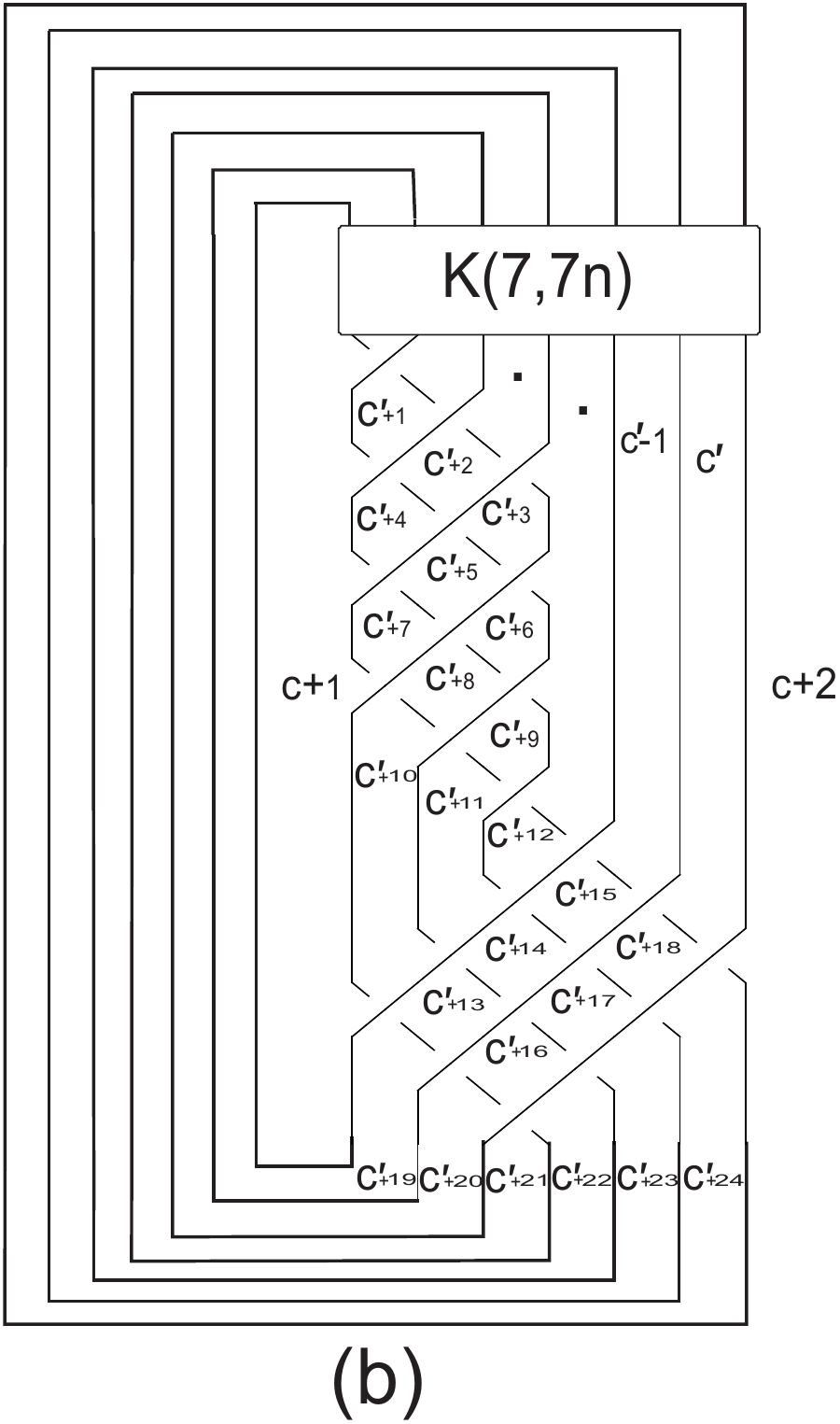}
\caption{Minimal diagrams for $(a)\ K(7,7)\ ;\ (b)\ K(7,7n+4)$}
\end{center}
\end{figure}
Let $X_i = \{2i(p-1), 2i(p-1)-2, \cdots, 2i(p-1)-2(i-1)\}$ be the set of $ i $ regions in $(\sigma_1\sigma_2\cdots\sigma_{p-1})^{q}$. We use $X_i^{'s}$ to identify the selected regions throughout the paper.
\begin{thm}\label{thm3} 
For a torus link $K(p,q)$, where $q = np\ or\ np+1$, we obtain the following: 
\begin{enumerate}
\item If $p$ is odd, then \[\displaystyle u_R(K(p,q))\leq \frac{n(p^2-1)}{8}.\]
\item If $p$ is even and $n$ is odd, then $K(p,np)$ is not proper and \[u_R(K(p,np+1))\leq \frac{np^2+2p}{8}.\]
\item If both $p$ and $n$ are even, then  \[u_R(K(p,q))\leq \frac{np^2}{8}.\]
\end{enumerate}
\end{thm}
\begin{pf} Consider $K(p,q)$ as the closure of $(\sigma_1\sigma_2\cdots\sigma_{p-1})^{q}$, as shown in Figure 1(a). 
Then,\\
\textbf{Case 1}. If $p$ is odd, then we select $\{(i-1)(p-1)p + \{X_1, X_2,\cdots, X_{\frac{p-1}{2}}\}\}_{i=1}^n$ regions so that after making region crossing change on these regions, we obtain
\begin{align}
(\mu_1\mu_2\cdots\mu_{p})^n &; if q = np, \notag \\
(\mu_1\mu_2\cdots\mu_{p})^n\mu_1 &; if q = np+1 \notag
\end{align}
where $\mu_i = \sigma_1\sigma_2\cdots\sigma_{p-i}\sigma_{p-i+1}^{-1}\sigma_{p-i+2}^{-1}\cdots\sigma_{p-1}^{-1}$. But 
$(\mu_1\mu_2\cdots\mu_{p})^n$ is a trivial $p$-braid by Theorem \ref{thm2}. Hence, the closure of the resultant braid is a trivial link.

Observe that the number of region crossing changes is equal to $n\sum(\frac{p-1}{2}) = \frac{n(p^2-1)}{8}$. So
\[u_R(K(p,np))\leq \frac{n(p^2-1)}{8}.\]
\textbf{Case 2}. If $p$ is even and $n$ is odd, then that $K(p,np)$ is not proper follows from Theorem \ref{thm1} and for $K(p,np+1)$ 
we select the regions $\{(i-1)(p-1)2p + \{X_1, X_2,\cdots, X_{\frac{p}{2}},((p-1)(2p-1)+1)-\{X_1, X_2,\cdots, X_{\frac{p-2}{2}}\}\}\}_{i=1}^{\frac{n-1}{2}}, (n-1)p(p-1)+\{X_1, X_2,\cdots, X_{\frac{p}{2}}\}$. Then, after making region crossing change on these regions, we obtain
\[(\mu_1\mu_2\cdots\mu_{p}\nu_1\nu_2\cdots\nu_{p})^\frac{n-1}{2}\mu_1\mu_2\cdots\mu_{p-1}\mu_p\mu_p\]
where $\mu_i = \sigma_1\sigma_2\cdots\sigma_{p-i}\sigma_{p-i+1}^{-1}\sigma_{p-i+2}^{-1}\cdots\sigma_{p-1}^{-1}$ and
$\nu_i = \sigma_1^{-1}\sigma_2^{-1} \cdots \sigma_{p-i}^{-1}\sigma_{p-i+1} \cdots\sigma_{p-1}$. Note that 
$\nu_i$ is a mirror image of $\mu_i$. Observe that $\mu_1\mu_2\cdots\mu_{p}$ and $\nu_1\nu_2\cdots\nu_{p}$ are trivial $p$-braids by Theorem \ref{thm2} and Remark \ref{coro2}. Hence, the closure of the resultant braid is a trivial link.

Observe that the number of region crossing changes is equal to $\frac{(n-1)}{2}(\sum\frac{p}{2} + \sum \frac{p-2}{2})+
\sum\frac{p}{2} = \frac{np^2+2p}{8}$. Thus, \[u_R(K(p,np+1))\leq \frac{np^2+2p}{8}.\]
\textbf{Case 3}. If both $p$ and $n$ are even, then we select $\{(i-1)(p-1)2p + \{X_1, X_2,\cdots, X_{\frac{p}{2}},((p-1)(2p-1)+1)-\{X_1, X_2,\cdots, X_{\frac{p-2}{2}}\}\}\}_{i=1}^{\frac{n}{2}}$ regions so that after making region crossing change on these regions, we obtain
\begin{align}
(\mu_1\mu_2\cdots\mu_{p}\nu_1\nu_2\cdots\nu_{p})^\frac{n}{2} &; if q = np, \notag \\
(\mu_1\mu_2\cdots\mu_{p}\nu_1\nu_2\cdots\nu_{p})^\frac{n}{2}\mu_1 &; if q = np+1 \notag
\end{align}
where $\mu_i = \sigma_1\sigma_2\cdots\sigma_{p-i}\sigma_{p-i+1}^{-1}\sigma_{p-i+2}^{-1}\cdots\sigma_{p-1}^{-1}$ and
$\nu_i = \sigma_1^{-1}\sigma_2^{-1} \cdots \sigma_{p-i}^{-1} \sigma_{p-i+1} \cdots\sigma_{p-1}$. So by Theorem~\ref{thm2} and Remark \ref{coro2}, the closure of the resultant braid is a trivial link.

Observe that the number of region crossing changes is equal to $\frac{n}{2}(\sum\frac{p}{2} + \sum \frac{p-2}{2}) = \frac{np^2}{8}$. Thus \\ $~~~~~~~~~~~~~~~~~~~~~~~~~~~~~~~~~~~~~~~~~~~~~~~~~~\displaystyle u_R(K(p,np))\leq \frac{np^2}{8}.$ \hfill $\square$
\end{pf} 

Observe that the unknotting number of $K(p,np)$ and $K(p,np+1)$ are identical. Also, the number of crossing changes due to the region crossing changes prescribed in Theorem \ref{thm3} for proper $K(p,np)$ and $K(p,np+1)$ torus links is equal to the unknotting number for $K(p,np)$.
\begin{rmk}\label{rmk2}
Consider a proper $K(2,q)$ torus knot/link then, the only possible minimal diagrams for $K(2,q)$ are the closures of $\sigma_1^q$, as shown in Figure 2. 
In any case, the region crossing change at $(q+1)^{th}$ or  $(q+2)^{th}$ region provides us the mirror image of $K(2,q)$. The region crossing change at any other region provides us $K(2,q-4)$. Hence, it is easy to observe that $u_R(K(2,q)) = \lfloor\frac{q+2}{4}\rfloor$.
\begin{figure}\label{fig2}
\begin{center}
\includegraphics[width=4cm,height=2.5cm]{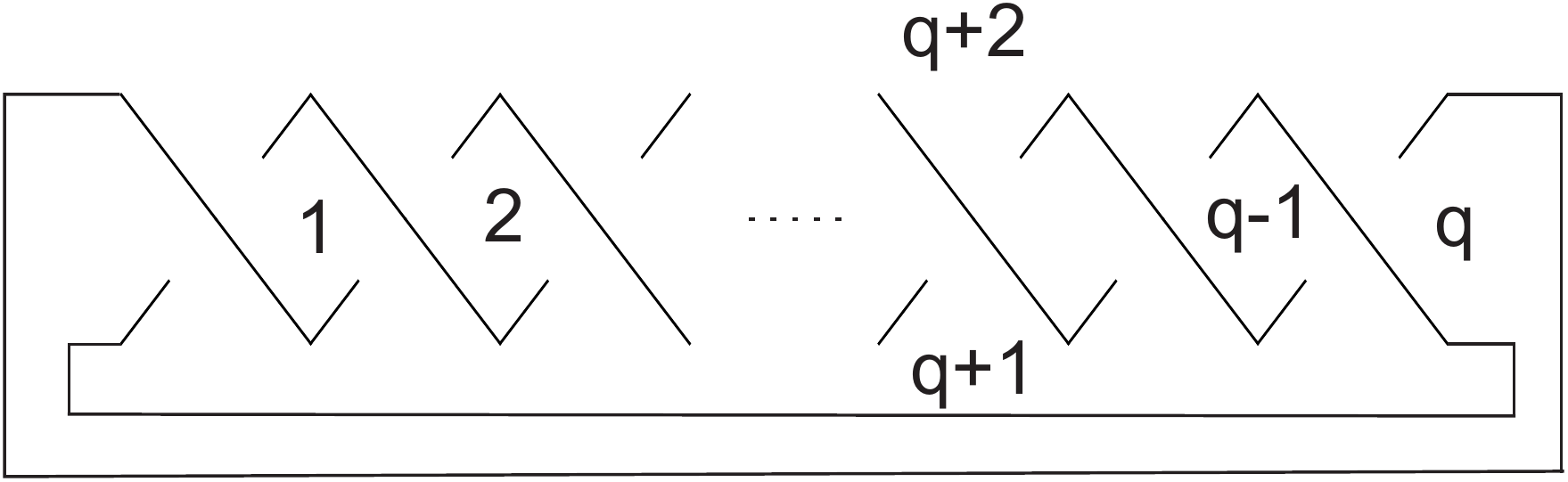}~~~~~
\includegraphics[width=4cm,height=2cm]{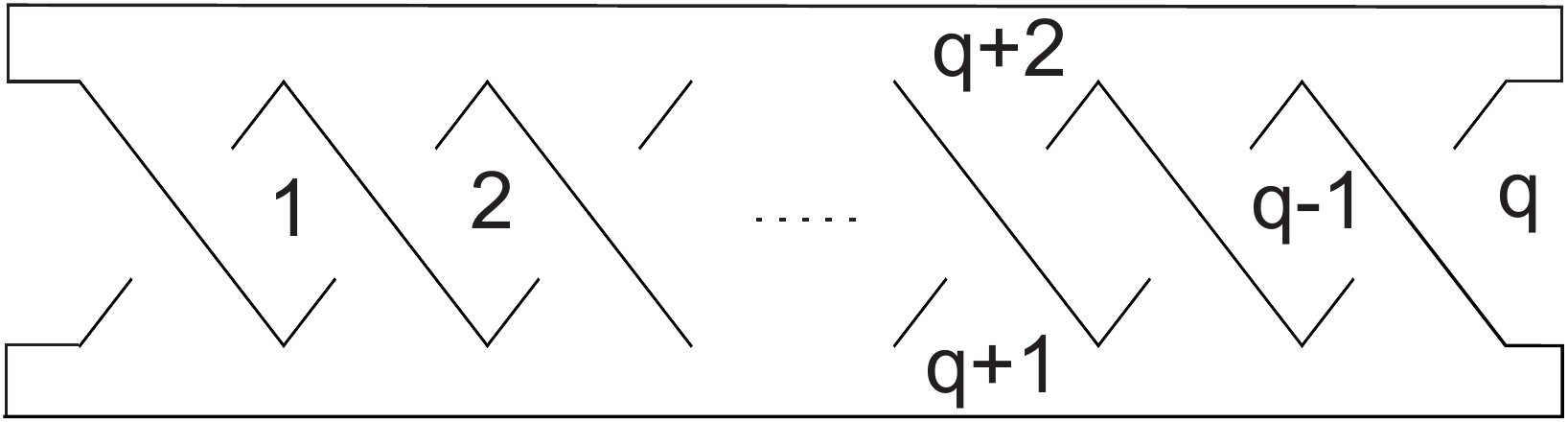}
\caption{Minimal diagrams for $ K(2,q)$}
\end{center}
\end{figure}
\end{rmk}
\begin{thm}\label{thm16} For torus links $K(p,np+a)$, where $a$ is odd and $p\equiv 0\ or\ \pm1\ (mod\ a)$, $n\geq1$. Then, we obtain the following:
\begin{enumerate}
\item If $p$ is odd, then 
\[ u_R(K(p,np+a))\leq \frac{n(p^2-1)}{8}+\left\lfloor\frac{p+1}{a}\right\rfloor\left(\frac{a^2-1}{8}\right).\]
\item If both $p$ and $n$ are even, then
\[u_R(K(p,np+a))\leq \frac{np^2}{8}+\left\lfloor\frac{p+1}{a}\right\rfloor\left(\frac{a^2-1}{8}\right).\]
\end{enumerate}
\end{thm}
\begin{pf} Consider $K(p,np+a)$ as the closure of $(\sigma_1\sigma_2\cdots\sigma_{p-1})^{np+a}$ as in Figure 1(a). 
Then,\\
\textbf{Case 1}. If $p$ is odd, then we select  $\{\{(i-1)(p-1)p + \{X_1, X_2,\cdots, X_{\frac{p-1}{2}}\}\}_{i=1}^n, \{(np(p-1)-ja)+\{X_1, X_2,\cdots, X_{\frac{a-1}{2}}\}\}_{j=0}^{\lfloor\frac{p+1}{a}\rfloor-1}\}$ regions so that after making region crossing changes at these regions, the closure of the resultant braid is a trivial link by
Theorem \ref{thm3}, Remark \ref{coro2} and Lemma \ref{thm7}. Observe that the number of the regions selected is $n\sum\frac{p-1}{2}+ \lfloor\frac{p+1}{a}\rfloor\sum\frac{a-1}{2} = \frac{n(p^2-1)}{8}+\lfloor\frac{p+1}{a}\rfloor\frac{a^2-1}{8}$. Thus, \[u_R(K(p,np+a))\leq \frac{n(p^2-1)}{8}+\left\lfloor\frac{p+1}{a}\right\rfloor\left(\frac{a^2-1}{8}\right).\]
\textbf{Case 2}. If both $p$ and $n$ are even, then we select
$\{\{(i-1)(p-1)2p + \{X_1, X_2,\cdots, X_{\frac{p}{2}},((p-1)(2p-1)+1)-\{X_1, X_2,\cdots, X_{\frac{p-2}{2}}\}\}\}_{i=1}^{\frac{n}{2}}, \{(np(p-1)-ja)+\{X_1, X_2,\cdots, X_\frac{a-1}{2}\}\}_{j=0}^{\lfloor\frac{p+1}{a}\rfloor-1}\}$ regions so that after making region crossing changes at these regions, the closure of the resultant braid is a trivial link by
Theorem \ref{thm3}, Remark \ref{coro2} and Lemma \ref{thm7}. Observe that the number of regions selected is $\frac{n}{2}(\sum\frac{p}{2}+ \sum\frac{p-2}{2})+\lfloor\frac{p+1}{a}\rfloor\sum\frac{a-1}{2} = \frac{np^2}{8}+\lfloor\frac{p+1}{a}\rfloor\frac{a^2-1}{8}$. Thus, \\
$\displaystyle ~~~~~~~~~~~~~~~~~~~~~~~~~~~~~~~~~~~u_R(K(p,np+a))\leq \frac{np^2}{8}+\left\lfloor\frac{p+1}{a}\right\rfloor\left(\frac{a^2-1}{8}\right).$ \hfill $\square $
\end{pf}
\begin{thm}\label{p2moda} For torus links $K(p,np+a)$, where $a$ is odd and $p\equiv 2\ (mod\ a)$, $n\geq1$, we obtain the following:
\begin{enumerate}
\item If $p$ is odd, then 
\[u_R(K(p,np+a))\leq \frac{n(p^2-1)}{8}+\frac{p-2}{a}\left(\frac{a^2-1}{8}\right)+\left\lfloor\frac{a+2}{4}\right\rfloor.\]
\item If both $p$ and $n$ are even, then
\[u_R(K(p,np+a))\leq \frac{np^2}{8}+\frac{p-2}{a}\left(\frac{a^2-1}{8}\right)+\left\lfloor\frac{a+2}{4}\right\rfloor.\]
\end{enumerate}
\end{thm}
\begin{pf} Consider $K(p,np+a)$ as the closure of $(\sigma_1\sigma_2\cdots\sigma_{p-1})^{np+a}$ as in Figure 1(a).  
 Then,\\
\textbf{Case 1}. If both $p$ and $a$ are odd, then $\exists\ m\in \mathbb{N} \cup \{0\}$ such that $a= 4m\pm 1$. We select  $\{\{(i-1)(p-1)p + \{X_1, X_2,\cdots, X_{\frac{p-1}{2}}\}\}_{i=1}^n, \{(np(p-1)-ja)+\{X_1, X_2,\cdots, X_{\frac{a-1}{2}}\}\}_{j=0}^{\frac{p-2}{a}-1}, \{(np+1)(p-1)+1+4k(p-1)\}_{k=0}^{m-1}\}$ regions so that after making region crossing changes at these regions, the closure of the resultant braid is a trivial link by
Theorem \ref{thm3} and Lemma \ref{thm7}. Observe that the number of selected regions is $n\sum\frac{p-1}{2}+ \frac{p-2}{a}\sum\frac{a-1}{2}+m = n\frac{p^2-1}{8}+\frac{p-2}{a}(\frac{a^2-1}{8}) +\lfloor\frac{a+2}{4}\rfloor$. Thus, \[u_R(K(p,np+a))\leq \frac{n(p^2-1)}{8}+\frac{p-2}{a}\left(\frac{a^2-1}{8}\right)+\left\lfloor\frac{a+2}{4}\right\rfloor.\]
\textbf{Case 2}. If both $p$ and $n$ are even and $a$ is odd, then $\exists\ m\in \mathbb{N} \cup \{0\}$ such that $a= 4m\pm 1$. We select $\{\{(i-1)(p-1)2p + \{X_1, X_2,\cdots, X_{\frac{p}{2}},((p-1)(2p-1)+1)-\{X_1, X_2,\cdots, X_{\frac{p-2}{2}}\}\}\}_{i=1}^{\frac{n}{2}}, \{(np(p-1)-ja)+\{X_1, X_2,\cdots, X_\frac{a-1}{2}\}\}_{j=0}^{\frac{p-2}{a}-1}, \{(np+1)(p-1)+1+4k(p-1)\}_{k=0}^{m-1}\}$ regions so that after making region crossing changes at these regions, the closure of the resultant braid is a trivial link by
Theorem \ref{thm3} and Lemma \ref{thm7}. Observe that the number of selected regions is $\frac{n}{2}(\sum\frac{p}{2}+ \sum\frac{p-2}{2})+\frac{p-2}{a}\sum\frac{a-1}{2}+m = n\frac{p^2}{8}+\frac{p-2}{a}(\frac{a^2-1}{8})+\lfloor\frac{a+2}{4}\rfloor$. Thus,\\ $\displaystyle ~~~~~~~~~~~~~~~~~~~~~~~~~~~~~~~~~~u_R(K(p,np+a))\leq n\frac{p^2}{8}+\frac{p-2}{a}\left(\frac{a^2-1}{8}\right)+\left\lfloor\frac{a+2}{4}\right\rfloor.$ \hfill $\square$
\end{pf}
%
\begin{thm}\label{thm17} For torus links $K(p,np+a)$ where $a$ is even, $n$ is odd and $p\equiv 0\ (mod\ a)$, we obtain the following:
\[u_R(K(p,np+a))\leq \frac{np^2+ap}{8}.\]
\end{thm}
\begin{pf}
Consider $K(p,np+a)$ as the closure of $(\sigma_1\sigma_2\cdots\sigma_{p-1})^{np+a}$ as in Figure 1(a). 
 Then, we select
$\{\{(i-1)(p-1)2p + \{X_1, X_2,\cdots, X_{\frac{p}{2}},((p-1)(2p-1)+1)-\{X_1, X_2,\cdots, X_{\frac{p-2}{2}}\}\}\}_{i=1}^{\frac{n-1}{2}}, (n-1)p(p-1)+\{X_1, X_2,\cdots, X_\frac{p}{2}, \{ja+(p+a-1)(p-1)+1-\{X_1, X_2,\cdots, X_\frac{a-2}{2}\}\}_{j=1}^{\frac{p}{a}-1}\}\}$ regions so that after making region crossing changes at these regions, the closure of the resultant braid is a trivial link by Theorem \ref{thm3} and Lemma \ref{thm10}. Observe that the number of regions selected is $\frac{n-1}{2}\left(\sum(\frac{p}{2})+ \sum(\frac{p-2}{2})\right)+\sum(\frac{p}{2})+\frac{p}{a}\sum(\frac{a-2}{2}) = \frac{(n-1)p^2}{8}+\frac{p(p+2)}{8}+\frac{p}{a}\left(\frac{a(a-2)}{8}\right)$. Thus,\\ $\displaystyle ~~~~~~~~~~~~~~~~~~~~~~~~~~~~~~~~~~~~~~~~~~~~~~u_R(K(p,np+a))\leq \frac{np^2+ap}{8}.$  \hfill $\square$
\end{pf}
\begin{thm}\label{thm5} For torus links $K(p,np-1)$ with $n\geq2$
\begin{enumerate}
\item If $p$ is odd, then \[u_R(K(p,np-1))\leq \frac{n(p^2-1)}{8}.\]
\item If $p$ is even and $n$ is odd, then\[u_R(K(p,np-1))\leq \frac{np^2-2p}{8}.\]
\item If both $p$ and $n$ are even, then  \[u_R(K(p,np-1))\leq \frac{np^2}{8}.\]
\end{enumerate}
\end{thm}
\begin{pf} Consider $K(p,np-1)$ as the closure of $(\sigma_1\sigma_2\cdots\sigma_{p-1})^{np-1}$ as in Figure 1(a). 
Then,\\
\textbf{Case 1}. If $p$ is odd, then 
we select the regions $\{(i-1)(p-1)p + \{X_1, X_2,\cdots, X_{\frac{p-1}{2}} \}\}_{i=1}^{n-1}, ((n-1)p-1)(p-1) + \{X_1, X_2,\cdots, X_{\frac{p-1}{2}}\}$ so that after making region crossing change on these regions, we obtain
\[(\mu_1\mu_2\cdots\mu_{p})^{n-1}\mu_2\mu_3\cdots\mu_{p}\]
where $\mu_i = \sigma_1\sigma_2\cdots\sigma_{p-i}\sigma_{p-i+1}^{-1}\sigma_{p-i+2}^{-1}\cdots\sigma_{p-1}^{-1}$. By Theorem \ref{thm2} and Remark \ref{coro2}, the closure of $(\mu_1\mu_2\cdots\mu_{p})^{n-1} \mu_2\mu_3\cdots\mu_{p}$ is a trivial knot.

Observe that number of crossing changes is $(n-1)\sum\frac{p-1}{2} + \sum\frac{p-1}{2} = \frac{n(p^2-1)}{8}$. Thus, \[u_R(K(p,np+1))\leq \frac{n(p^2-1)}{8}\]
\textbf{Case 2}. If $p$ is even and $n$ is odd, then we write $K(p,np-1) = K(p, (n-1)p + (p-1))$, where both $p$ and $n-1$ are even. This result follows from Theorem \ref{thm16}.\\
%
\textbf{Case 3}. If both $p$ and $n$ are even, we select the regions same as we selected for $K(p,np)$ in Theorem \ref{thm3}. After making region crossing changes on selected regions, we are left with
 \[(\mu_1\mu_2\cdots\mu_{p}\nu_1\nu_2\cdots\nu_{p})^\frac{n-2}{2}\mu_1\mu_2\cdots\mu_{p}\nu_1\nu_2\cdots\nu_{p-1}.\]
By Theorem \ref{thm2},  Remark \ref{coro2} and Lemma \ref{thm7}, its closure is unknot. Observe that the number of region crossing changes is $n\left(\frac{p^2}{8}\right)$. Thus, \\ $\displaystyle ~~~~~~~~~~~~~~~~~~~~~~~~~~~~~~~~~~~~~~~~~~~~~~~~u_R(K(p,np))\leq \frac{np^2}{8}.$  \hfill $\square$
\end{pf}
\begin{thm}\label{thm8} For torus links $K(p,np+2)$ with $n\geq1$
\begin{enumerate}
\item If $p$ is odd, then
\[u_R(K(p,np+2))\leq \frac{n(p^2-1)}{8}+\left\lfloor\frac{p+1}{4}\right\rfloor.\]
\item If $p$ is even and $n$ is odd, then \[u_R(K(p,np+2))\leq \frac{np^2+2p}{8}.\]
\item If both $p$ and $n$ are even, then
\begin{enumerate}
\item if $p = 4m + 2$ for some $m$, then $K(p,np+2)$ is not proper.
\item if $p = 4m$ for some $m$, then\[u_R(K(p,np+2))\leq \frac{np^2+2p}{8}.\]
\end{enumerate}  
\end{enumerate}
\end{thm}

\begin{pf} Consider $K(p,np+2)$ as the closure of $(\sigma_1\sigma_2\cdots\sigma_{p-1})^{np+2}$ as in Figure 1(a). 
 Then,\\
\textbf{Case 1}. If $p$ is odd then $\exists\ m\in \mathbb{N} \cup \{0\}$ such that $p = 4m\pm1$. We select $\{(i-1)(p-1)p + \{X_1, X_2,\cdots, X_{\frac{p-1}{2}}\}\}_{i=1}^n, \{(np+1)(p-1)-4j\}_{j=0}^{m-1}$ regions so that after making region crossing change on these regions, we obtain by Theorem \ref{thm2} and Lemma \ref{prop1} that\\
\[\sigma_1\sigma_2\sigma_3\sigma_4^{-1}\sigma_1\sigma_2\sigma_3^{-1}\sigma_4^{-1};\ if\ p=4m+1,\]
 \[\sigma_1\sigma_2^{-1}\sigma_1^{-1}\sigma_2^{-1};\ if\ p=4m-1.\]
We observe that the closure of these braids are trivial links.

Observe that the number of region crossing changes is ${n}\sum(\frac{p-1}{2}) + m = n(\frac{p^2-1}{8}) + \lfloor\frac{p+1}{4}\rfloor$. Thus, \[u_R(K(p,np+2))\leq \frac{n(p^2-1)}{8}+\left\lfloor\frac{p+1}{4}\right\rfloor.\]
\textbf{Case 2}. If $p$ is even and $n$ is odd, then by Theorem \ref{thm17}:
%
\[u_R(K(p,np+2))\leq \frac{np^2+2p}{8}.\]
\textbf{Case 3}. If both $p$ and $n$ are even and $p\equiv 0\ (mod\ 4)$, then we select the regions $\{(i-1)(p-1)2p + \{X_1, X_2,\cdots, X_{\frac{p}{2}},((p-1)(2p-1)+1)-\{X_1, X_2,\cdots, X_{\frac{p-2}{2}}\}\}\}_{i=1}^{\frac{n}{2}}, \{(np+1)(p-1)-4j\}_{j=0}^{\frac{p}{4}-1}$ so that by Theorem \ref{thm2} and Lemma \ref{prop1}, after making region crossing change on these regions, we obtain
\[\sigma_1\sigma_2\sigma_3^{-1}\sigma_1\sigma_2^{-1}\sigma_3^{-1}.\]
We observe that the closure of this braid is a trivial link. Also, the number of region crossing changes is $\frac{n}{2}\left(\sum\frac{p}{2} + \sum \frac{p-2}{2}\right) + \frac{p}{4} = n(\frac{p^2}{8}) + \frac{p}{4}$. Thus, \[u_R(K(p,np+2))\leq \frac{np^2+2p}{8}.\]
But if $p\equiv 2\ (mod\ 4)$, then it follows from Theorem \ref{thm1} that $K(p,np+2)$ is not proper. \hfill $\square$ 
\end{pf}
\begin{thm}\label{p-2moda} Consider torus links $K(p,np+a)$ where $a$ is odd and $p\equiv -2\ (mod\ a)$, $n\geq1$. Then we obtain the following: 
\begin{enumerate}
\item If $p$ is odd, then 
\[u_R(K(p,np+a))\leq \frac{n(p^2-1)}{8}+\left(\frac{p+2}{a}-1\right)\left(\frac{a^2-1}{8}\right)+\frac{(a-2)^2-1}{8}+\left\lfloor\frac{a}{4}\right\rfloor.\]
\item If both $p$ and $n$ are even, then
\[u_R(K(p,np+a))\leq \frac{np^2}{8}+\left(\frac{p+2}{a}-1\right)\left(\frac{a^2-1}{8}\right)+\frac{(a-2)^2-1}{8}+\left\lfloor\frac{a}{4}\right\rfloor.\]
\end{enumerate}
\end{thm}
\begin{pf} Consider $K(p,np+a)$ as the closure of $(\sigma_1\sigma_2\cdots\sigma_{p-1})^{np+a}$ as in Figure 1(a). 
 Then,\\
\textbf{Case 1}. If $p$ is odd, then we select  $\{\{(i-1)(p-1)p + \{X_1, X_2,\cdots, X_{\frac{p-1}{2}}\}\}_{i=1}^n, \{(np(p-1)-ja)+\{X_1, X_2,\cdots, X_{\frac{a-1}{2}}\}\}_{j=0}^{\frac{p+2}{a}-2}, \{(np(p-1)-(\frac{p+2}{a}-1)a)+\{X_1, X_2,\cdots, X_{\frac{a-3}{2}}\}\}, \{(np+a-2)(p-1)+(a-3)-4l\}_{l=0}^{\lfloor\frac{a}{4}\rfloor}$ regions so that after making region crossing changes at these regions, the closure of the resultant braid is a trivial link by
Theorem \ref{thm3}, Lemma \ref{thm7} and Theorem \ref{thm8}. Observe that the number of selected regions  is $n\sum(\frac{p-1}{2})+ \left(\frac{p+2}{a}-1\right)\left(\sum(\frac{a-1}{2})\right)+\sum(\frac{a-3}{2})+\lfloor\frac{a}{4}\rfloor = \frac{n(p^2-1)}{8}+\left(\frac{p+2}{a}-1\right)\left(\frac{a^2-1}{8}\right)+\frac{(a-2)^2-1}{8} +\lfloor\frac{a}{4}\rfloor$. Thus, \[u_R(K(p,np+a))\leq \frac{n(p^2-1)}{8}+\left(\frac{p+2}{a}-1\right)\left(\frac{a^2-1}{8}\right)+\frac{(a-2)^2-1}{8}+\left\lfloor\frac{a}{4}\right\rfloor.\]
\textbf{Case 2}. If both $p$ and $n$ are even, then we select
$\{\{(i-1)(p-1)2p + \{X_1, X_2,\cdots, X_{\frac{p}{2}},((p-1)(2p-1)+1)-\{X_1, X_2,\cdots, X_{\frac{p-2}{2}}\}\}\}_{i=1}^{\frac{n}{2}},\{(np(p-1)-ja)+\{X_1, X_2,\cdots, X_{\frac{a-1}{2}}\}\}_{j=0}^{\frac{p+2}{a}-2}, \{(np(p-1)-(\frac{p+2}{a}-1)a)+\{X_1, X_2,\cdots, X_{\frac{a-3}{2}}\}\}, \{(np+a-2)(p-1)+(a-3)-4l\}_{l=0}^{\lfloor\frac{a}{4}\rfloor}$ regions so that after making region crossing changes at these regions, the closure of the resultant braid is a trivial link by
Theorem \ref{thm3}, Lemma \ref{thm7} and Theorem \ref{thm8}. Observe that the number of selected regions  is $\frac{n}{2}\left(\sum\frac{p}{2}+ \sum\frac{p-2}{2}\right)+\left(\frac{p+2}{a}-1\right)\left(\sum(\frac{a-1}{2})\right)+\sum(\frac{a-3}{2})+\lfloor\frac{a}{4}\rfloor = \frac{np^2}{8}+\left(\frac{p+2}{a}-1\right)\left(\frac{a^2-1}{8}\right)+\frac{(a-2)^2-1}{8} +\lfloor\frac{a}{4}\rfloor$. Thus,\\ $~~~~~~~~~~~~~~~~~~~~~~~\displaystyle  u_R(K(p,np+a))\leq \frac{np^2}{8}+\left(\frac{p+2}{a}-1\right)\left(\frac{a^2-1}{8}\right)+\frac{(a-2)^2-1}{8}+\left\lfloor\frac{a}{4}\right\rfloor.$ \hfill $\square $
\end{pf}
\begin{thm}\label{thm9} For torus links $K(p,np-2)$ with $n\geq2$, we obtain the following
\begin{enumerate}
\item If $p$ is odd, then
\[u_R(K(p,np-2))\leq \frac{n(p^2-1)}{8}-\frac{p-1}{2}+\left\lfloor\frac{p}{4}\right\rfloor.\]
\item If $p$ is even and $n$ is odd, then\[u_R(K(p,np-2))\leq \frac{np^2-2p}{8}.\]
\item If both $p$ and $n$ are even, then
\begin{enumerate}
\item if $p \equiv 2\ (mod\ 4)$, then $K(p,np-2)$ is not proper.
\item if $p \equiv 0\ (mod\ 4)$, then\[u_R(K(p,np-2))\leq \frac{np^2-2p}{8}.\]
\end{enumerate}  
\end{enumerate}
\end{thm}
\begin{pf} Consider $K(p,np-2)$ as the closure of $(\sigma_1\sigma_2\cdots\sigma_{p-1})^{np-2}$ as in Figure 1(a). 
Then\\
\textbf{Case 1}. If $p$ is odd, then result hold by Theorem \ref{p2moda}.\\
%
\textbf{Case 2}. If $p$ is even and $n$ is odd, then select the regions $\{(i-1)(p-1)2p + \{X_1, X_2,\cdots, X_{\frac{p}{2}},((p-1)(2p-1)+1)-\{X_1, X_2,\cdots, X_{\frac{p-2}{2}}\}\}\}_{i=1}^{\frac{n-1}{2}}, ((n-1)p-1)(p-1)+\{X_1, X_2,\cdots, X_{\frac{p-2}{2}}\}$. Then after making region crossing change on these regions, we obtain
\[(\mu_1\mu_2\cdots\mu_{p}\nu_1\nu_2\cdots\nu_{p})^\frac{n-1}{2}\mu_2\mu_3\cdots\mu_{p-1}\]
where $\mu_i = \sigma_1\sigma_2\cdots\sigma_{p-i}\sigma_{p-i+1}^{-1}\sigma_{p-i+2}^{-1}\cdots\sigma_{p-1}^{-1}$ and
$\nu_i = \sigma_1^{-1}\sigma_2^{-1} \cdots \sigma_{p-i}^{-1} \sigma_{p-i+1} \cdots\sigma_{p-1}$. By Theorem \ref{thm2}, Remark \ref{coro2} and Lemma \ref{prop2}, the closure of resultant braid is a trivial link.

Observe that the number of region crossing changes is equal to $\frac{(n-1)}{2}(\sum\frac{p}{2} + \sum \frac{p-2}{2})+\sum(\frac{p-2}{2}) = \frac{np^2-2p}{8}$. Thus, \[u_R(K(p,np+1))\leq \frac{np^2-2p}{8}.\]
\textbf{Case 3}. If both $p$ and $n$ are even, then
\begin{enumerate}
\item if $p=4m+2$ for some $m$, then it is clear from Theorem \ref{thm1} that $K(p,np-2)$ is not proper.
\item if $p=4m$ for some $m$, then we select the regions $\{(i-1)(p-1)2p + \{X_1, X_2,\cdots, X_{\frac{p}{2}},((p-1)(2p-1)+1)-\{X_1, X_2,\cdots, X_{\frac{p-2}{2}}\}\}\}_{i=1}^{\frac{n-2}{2}}, (n-2)p(p-1)+\{X_1, X_2,\cdots, X_{\frac{p}{2}}, (2p-3)(p-1)+2-\{X_1, X_2, \cdots, X_{\frac{p-4}{2}}\}, \{( p + 2+4i)(p-1)+1\}_{i=0}^{m-2}\}$. After making region crossing changes on selected regions, we are left with
 \[(\mu_1\mu_2\cdots\mu_{p}\nu_1\nu_2\cdots\nu_{p})^\frac{n-2}{2}\mu_1\mu_2\cdots\mu_{p}\beta\]
 where $\beta=\beta_1\beta_2\cdots\beta_{p-2}$ and $\beta_i$ is defined as:  $\beta_1=\sigma_1^{-1}\sigma_2^{-1}\cdots\sigma_{p-1}^{-1}$, $\beta_{p-2}=\sigma_1\sigma_2\cdots\sigma_{p-1}$ 
and for the remaining $i^{'s}$,\\ if $i\equiv 1\ (mod\ 4)$, 
$\beta_i = \sigma_1\sigma_2^{-1}\sigma_3^{-1}\sigma_4^{-1}\cdots\sigma_{p-i}^{-1}\sigma_{p-i+1}\cdots\sigma_{p-1}$,\\
if $i\equiv 2\ (mod\ 4)$, 
$\beta_i = \sigma_1\sigma_2\sigma_3^{-1}\sigma_4^{-1}\cdots\sigma_{p-1}^{-1}\sigma_{p-i+1}\cdots\sigma_{p-1}$ and\\
if $i\equiv 0\ or\ 3\ (mod\ 4)$, 
$\beta_i = \sigma_1^{-1}\sigma_2\sigma_3^{-1}\sigma_4^{-1}\cdots\sigma_{p-i}^{-1}\sigma_{p-i+1}\cdots\sigma_{p-1}$.\\
By Theorem \ref{thm2},  Remark \ref{coro2} and Lemma \ref{thm10}, the closure of\\
\[(\mu_1\mu_2\cdots\mu_{p}\nu_1\nu_2\cdots\nu_{p})^\frac{n-2}{2}\mu_1\mu_2\cdots\mu_{p}\beta\sim
\underbrace{\sigma_1^{-1}}_1\underbrace{\sigma_1}_2\underbrace{\sigma_1^{-1}}_3\underbrace{\sigma_1^{-1}}_4
\underbrace{\sigma_1}_5\cdots\underbrace{\sigma_1^{-1}}_{p-2}\]
is a trivial link. 

Observe that the number of region crossing changes is $\frac{n-2}{2}\left(\sum\frac{p}{2} + \sum \frac{p-2}{2}\right)
+ \sum(\frac{p}{2}) + \sum (\frac{p-4}{2}) + \frac{p}{4}-1 = \frac{np^2-2p}{8}$. Thus, \\ $\displaystyle ~~~~~~~~~~~~~~~~~~~~~~~~~~~~~~~~~~~~~~~~~~~~~~u_R(K(p,np))\leq \frac{np^2-2p}{8}.$ \hfill $\square$
\end{enumerate} 
\end{pf}
\begin{thm}\label{thm11} For torus links $K(p,np+3)$ with $n\geq1$
\begin{enumerate}
\item If $p$ is odd, then
\[u_R(K(p,np+3))\leq \frac{n(p^2-1)}{8}+\left\lfloor\frac{p+1}{3}\right\rfloor.\]
\item If $p$ is even and $n$ is odd, then 
\[u_R(K(p,np+3))\leq \frac{np^2+2p}{8}+\left\lfloor\frac{p+2}{6}\right\rfloor.\]
\item If both $p$ and $n$ are even, then
\[u_R(K(p,np+3))\leq \frac{np^2}{8}+\left\lfloor\frac{p+1}{3}\right\rfloor.\]
\end{enumerate}  
\end{thm}
\begin{pf}  Consider $K(p,np+3)$ as the closure of $(\sigma_1\sigma_2\cdots\sigma_{p-1})^{np+3}$ as in Figure 1(a). 
Then\\
\textbf{Case 1}. If $p$ is odd, then the result follows from Theorem \ref{thm16}.\\
%
\textbf{Case 2}. If $p$ is even and $n$ is odd, then $p+2 = 6m+a$ for some $m\geq0$ and $0\leq a<6$ and we select the regions  $\{(i-1)(p-1)2p + \{X_1, X_2,\cdots, X_{\frac{p}{2}},((p-1)(2p-1)+1)-\{X_1, X_2,\cdots, X_{\frac{p-2}{2}}\}\}\}_{i=1}^{\frac{n-1}{2}}, (n-1)p(p-1)+\{X_1, X_2,\cdots, X_{\frac{p}{2}}\}, \{(np+2)(p-1)-2-6i\}_{i=0}^{m-1}$. 
Then after making region crossing change on these regions, we obtain
\[(\mu_1\mu_2\cdots\mu_{p}\nu_1\nu_2\cdots\nu_{p})^\frac{n-1}{2}\mu_1\mu_2\cdots\mu_{p-1}\mu_p\beta,\]
where $\mu_i = \sigma_1\sigma_2\cdots\sigma_{p-i}\sigma_{p-i+1}^{-1}\sigma_{p-i+2}^{-1}\cdots\sigma_{p-1}^{-1}$, 
$\nu_i = \sigma_1^{-1}\sigma_2^{-1} \cdots \sigma_{p-i}^{-1} \sigma_{p-i+1}\cdots\sigma_{p-1}$ and $\beta=\sigma_1^{-1}\sigma_2^{-1}\\ \cdots\sigma_{p-1}^{-1}\sigma_1^{g_{2,1}}\sigma_2^{g_{2,2}}\cdots
\sigma_{p-1}^{g_{2,p-1}}\sigma_1^{g_{3,1}}\sigma_2^{g_{3,2}}\cdots
\sigma_{p-1}^{g_{3,p-1}}$ with 
$g_{2,j} = -1$ for $j = \{\{(p-2),(p-3)\}-6i\}_{i=0}^{m-1}$ and
$g_{3,j} = -1$ for $j = \{\{(p-3),(p-4)\}-6i\}_{i=0}^{m-1}$.\\
By Theorem \ref{thm2} and Lemma \ref{thm12}, the closure of the resultant braid is a trivial link.

Observe that the number of region crossing changes is equal to $\frac{(n-1)}{2}(\sum\frac{p}{2} + \sum \frac{p-2}{2})+ \sum(\frac{p}{2}) + m= \frac{np^2+2p}{8} +\lfloor\frac{p+2}{6}\rfloor$. Thus, \[u_R(K(p,np+3))\leq \frac{n(p^2+2p)}{8}+\left\lfloor\frac{p+2}{6}\right\rfloor.\]
\textbf{Case 3}. If both $p$ and $n$ are even, then the result follows from Theorem \ref{thm16}. \hfill $\square$
%
\end{pf} 
\begin{lem}\label{thm+4}
The closure of the braid\\
$\underbrace{\sigma_i^{-1}\sigma_{i-1}^{-1}\sigma_{i-2}^{-1}\sigma_{i-3}^{-1}}_1\underbrace{\sigma_{i+1}^{-1}\sigma_i^{-1}
\sigma_{i-1}^{-1}\sigma_{i-2}^{-1}}_2\underbrace{\sigma_{i+2}\sigma_{i+1}\sigma_{i}\sigma_{i-1}}_3\underbrace{
\sigma_{i+3}^{-1}\sigma_{i+2}^{-1}
\sigma_{i+1}\sigma_{i}}_4\underbrace{\sigma_{i+4}^{-1}\sigma_{i+3}^{-1}\sigma_{i+2}\sigma_{i+1}}_5  \underbrace{\sigma_{i+5}\sigma_{i+4}\sigma_{i+3}\sigma_{i+2}}_6\\ \underbrace{\sigma_{i+6}^{-1}\sigma_{i+5}^{-1}
\sigma_{i+4}\sigma_{i+3}}_7\underbrace{\sigma_{i+7}\sigma_{i+6}\sigma_{i+5}\sigma_{i+4}}_8$ is trivial link.
\end{lem}
\begin{pf} The proof follows elimentary braid relations.
\end{pf} 
\begin{thm}\label{thm15} For torus links $K(p,np+4)$ with $n\geq1$
\begin{enumerate}
\item If $p$ is odd, then
\begin{enumerate}
\item if $p\equiv 1\ or\ 3\ (mod\ 8)$, then 
\[u_R(K(p,np+4))\leq \frac{n(p^2-1)}{8}+\left\lfloor\frac{p}{2}\right\rfloor.\]
\item if $p\equiv 5\ or\ 7\ (mod\ 8)$, then
\[u_R(K(p,np+4))\leq \frac{n(p^2-1)}{8}+\left\lceil\frac{p}{2}\right\rceil.\]
\end{enumerate} 
\item If both $p$ and $n$ are even, then
\begin{enumerate}
\item if $p\equiv 0, 2\ or\ 6\ (mod\ 8)$, then 
\[u_R(K(p,np+4))\leq \frac{np^2}{8}+\frac{p}{2}.\]
\item if $p\equiv 4\ (mod\ 8)$, then $K(p,np)$ is not proper.
\end{enumerate}
\item If $p$ is even and $n$ is odd, then
\begin{enumerate}
\item if $p\equiv 0\ (mod\ 4)$, then
\[u_R(K(p,np+4))\leq \frac{np^2}{8}+ \frac{p}{2}.\]
\item if $p\equiv 2\ (mod\ 4)$, then $K(p,np+4)$ is not proper.
\end{enumerate}
\end{enumerate} 
\end{thm}
\begin{pf}
Consider $K(p,np+4)$ as the closure of $(\sigma_1\sigma_2\cdots\sigma_{p-1})^{np}(\sigma_1\sigma_2\sigma_{3})^{4}\sigma_4\sigma_{3}\cdots\sigma_1\sigma_{5}\sigma_{4}\cdots \sigma_2 \cdots\sigma_{p-1}\\ \sigma_{p-2}\cdots\sigma_{p-4}$ as in Figure 1(b). 
Also, if  $p\equiv a\ (mod\ 8)$, then $\exists\ m\in  \mathbb{N} \cup \{0\}$ such that $p= mp+a$. Then\\
\textbf{Case 1}. If $p$ is odd, then
\begin{enumerate}
\item if $p\equiv 1\ (mod\ 8)$, then we select $ \{\{(i-1)(p-1)p + \{X_1, X_2,\cdots, X_{\frac{p-1}{2}}\}\}_{i=1}^n,  np(p-1) + \{6,10,12,18, \{24i+\{1,3,12,18\}\}_{i=1}^{m-1}\}\}$ regions.
\item if $p\equiv 3\ (mod\ 8)$, then we select $ \{\{(i-1)(p-1)p + \{X_1, X_2,\cdots, X_{\frac{p-1}{2}}\}\}_{i=1}^n,  np(p-1) + \{6,10,12,18, \{24i+\{1,3,12,18\}\}_{i=1}^{m-1}, 24m+\{3\}\}\}$ regions.
\item if $p\equiv 5\ (mod\ 8)$ and $m\neq 0$, then we select $ \{\{(i-1)(p-1)p + \{X_1, X_2,\cdots, X_{\frac{p-1}{2}}\}\}_{i=1}^n, np(p-1) + \{6,10,12,18, \{24i+\{1,3,12,18\}\}_{i=1}^{m-1}, 24m+\{1,3,12\}\}\}$ regions and if $m=0$, then we select $\{\{(i-1)(p-1)p + \{X_1, X_2,\cdots, X_{\frac{p-1}{2}}\}\}_{i=1}^n, np(p-1) +\{6,10,12\}\}$ regions.
\item if $p\equiv 7\ (mod\ 8)$ and $m\neq 0$, then we select $ \{\{(i-1)(p-1)p + \{X_1, X_2,\cdots,  X_{\frac{p-1}{2}}\}\}_{i=1}^n, np(p-1) + \{6,10,12,18, \{24i+\{1,3,12,18\}\}_{i=1}^{m}\}\}$ regions and if $m=0$, then we select $\{\{(i-1)(p-1)p + \{X_1, X_2,\cdots, X_{\frac{p-1}{2}}\}\}_{i=1}^n, np(p-1) +\{6,10,12,18\}\}$ regions.
\end{enumerate}
\textbf{Case 2}. If both $p$ and $n$ are even, then 
\begin{enumerate}
\item if $p\equiv 0\ (mod\ 8)$, then we select $ \{\{(i-1)(p-1)2p + \{X_1, X_2,\cdots, X_{\frac{p}{2}},((p-1)(2p-1)+1)-\{X_1, X_2,\cdots, X_{\frac{p-2}{2}}\}\}\}_{i=1}^{\frac{n}{2}}, np(p-1) + \{6,10,12,18, \{24i+\{1,3,12,18\}\}_{i=1}^{m-1}\}\}$ regions.
\item if $p\equiv 2\ (mod\ 8)$, then we select $ \{\{(i-1)(p-1)2p + \{X_1, X_2,\cdots, X_{\frac{p}{2}},((p-1)(2p-1)+1)-\{X_1, X_2,\cdots, X_{\frac{p-2}{2}}\}\}\}_{i=1}^{\frac{n}{2}}, np(p-1) + \{6,10,12,18, \{24i+\{1,3,12,18\}\}_{i=1}^{m-1}, 24m+\{3\}\}\}$ regions.
\item if $p\equiv 6\ (mod\ 8)$ and $m\neq 0$, then we select $ \{\{(i-1)(p-1)2p + \{X_1, X_2,\cdots, X_{\frac{p}{2}},((p-1)(2p-1)+1)-\{X_1, X_2,\cdots, X_{\frac{p-2}{2}}\}\}\}_{i=1}^{\frac{n}{2}}, np(p-1) + \{6,10,12,18, \{24i+\{1,3,12,18\}\}_{i=1}^{m-1}, 24m+\{3\}\}\}$ regions and if $m=0$, then we select $ \{\{(i-1)(p-1)2p + \{X_1, X_2,\cdots, X_{\frac{p}{2}},((p-1)(2p-1)+1)-\{X_1, X_2,\cdots, X_{\frac{p-2}{2}}\}\}\}_{i=1}^{\frac{n}{2}}, np(p-1) + \{6,10,12\}\}$ regions.
\item if $p\equiv 4\ (mod\ 8)$, then it is not proper follows from Theorem \ref{thm1}.
\end{enumerate}
In all the above cases, after making region crossing changes the closure of the resultant braid is trivial link by Theorem \ref{thm3} and Lemma \ref{thm+4}.\\
\textbf{Case 3}. If $p$ is even and $n$ is odd, then
\begin{enumerate}
\item if $p\equiv 0\ (mod\ 4)$, then the result follows from Theorem \ref{thm17}.
\item if $p\equiv 2\ (mod) 4$, then the result follows from Theorem \ref{thm1}.\hfill $\square$
\end{enumerate} 
\end{pf}
\begin{thm}\label{thm15} For torus links $K(p,np+5)$ with $n\geq1$
\begin{enumerate}
\item If $p$ is odd, and
\begin{enumerate}
\item if $p\equiv 0\ or\ \pm 1\ (mod\ 5)$, then 
\[u_R(K(p,np+5))\leq \frac{n(p^2-1)}{8}+3\left\lfloor\frac{p+1}{5}\right\rfloor.\]
\item if $p\equiv 2\ (mod\ 5)$, then
\[u_R(K(p,np+5))\leq \frac{n(p^2-1)}{8}+\frac{3p-1}{5}.\]
\item if $p\equiv 3\ (mod\ 5)$, then
\[u_R(K(p,np+5))\leq \frac{n(p^2-1)}{8}+\frac{3p+1}{5}.\]
\end{enumerate} 
\item If both $p$ and $n$ are even, and
\begin{enumerate}
\item if $p\equiv 0\ or\ \pm 1\ (mod\ 5)$, then 
\[u_R(K(p,np+5))\leq \frac{np^2}{8}+3\left\lfloor\frac{p+1}{5}\right\rfloor.\]
\item if $p\equiv 2\ (mod\ 5)$, then
\[u_R(K(p,np+5))\leq \frac{np^2}{8}+\frac{3p-1}{5}.\]
\item if $p\equiv 3\ (mod\ 5)$, then
\[u_R(K(p,np+5))\leq \frac{np^2}{8}+\frac{3p+1}{5}.\]
\end{enumerate}
\end{enumerate} 
\end{thm}
\begin{pf} The proof directly follows from Theorem \ref{thm16}, Theorem \ref{p2moda} and Theorem \ref{p-2moda}.\hfill $\square$
\end{pf} 

\section{Conclusion}
It would be interesting to extend the results obtained here to find exact region unknotting number for torus knots/links. This seems to be a hard problem. In \cite{ayaka}, Ayaka showed that $u_R(K)\leq c(K)/2+1$, where $C(K)$ is the number of crossings of the knot K. Here, the sharp bounds provided for torus knots/links are much smaller than the upper bound $c(K)/2+1$. In particular, whatever be the case we considered in this paper, the upper bound we obtain, for $u_R(K)$, lies between  ${c}/{8}$ and ${c}/{5}$ and as $p \rightarrow \infty$, the sharp upper bound provided tends to $c/8$. In case of proper $K(p,np)$ and $K(p,np+1)$, we conjecture that the sharp upper bound provided in the paper is equal to the region unknotting number for them.\\

\noindent\textbf{Acknowledgements}\\

The first author thanks CSIR New Delhi and IIT Ropar for providing financial assistance and research facilities.


\begin{thebibliography}{99} 
%
\bibitem{ps1}
V. Siwach and P. Madeti, A method for unknotting torus knots. math.GT/1207.4918v1, 2012.

\bibitem{ayaka}
Ayaka Shimizu, Region crossing change is an unknotting operation. math.GT/1011.6304v2, 2010.


\bibitem{ayaka2}
Ayaka Shimizu, On region unknotting numbers. RIMS Kokyuroku 1766, 15-22, 2011.  

\bibitem{hash}
H. Murakami, Some metrics on classical knots. Math. Ann. 270, 35-45, 1985.

\bibitem{km-2}
P. Kronheimer and T. Mrowka, Gauge theory for embedded surfaces, I,
Topology, 32 , 773-826. 1993.


\bibitem{km-3}
P. Kronheimer and T. Mrowka, Gauge theory for embedded surfaces, II,
Topology, 34 , 37-97, 1995.

\bibitem{delta}
H. Murakami, Y. Nakanishi, On a certain move generating link-homology. Math. Ann. 284, 75-89, 1989.


\bibitem{incidence}
Cheng Zhiyun, Gao Hongzhu, On region crossing change and incidence matrix, Science China Mathematics. Vol. 55 No. 7. 1487–1495, 2012.


\bibitem{3-gon}
Y. Nakanishi, Replacements in the Conway third identity, Tokyo J. Math. 14, 197-203, 1991.

\bibitem{n-gon}
Haruko Aida, Unknotting operation for Polygonal type, Tokyo J. Math. Vol. 15, No. 1, 111-121, 1992.

\bibitem{H(n)}
J. Hoste, Y. Nakanishi and K. Taniyama, Unknotting operations involving trivial tangles, Osaka J.
Math. 27 , 555-566, 1990.

\bibitem{proper}
Cheng Zhiyun, When is region crossing change an unknotting operation? math.GT/1201.1735v1, 2012.

\end{thebibliography}
\end{document}